\documentclass[11pt]{article}
\usepackage{theorem}
\usepackage[english]{babel}
\usepackage{a4}
\usepackage{epsfig}
\usepackage{amsfonts}
\usepackage{latexsym}
\usepackage{html}
\usepackage{amsmath}

\newtheorem{theorem}{Theorem}[section]
\newtheorem{lemma}[theorem]{Lemma}
\newtheorem{proposition}[theorem]{Proposition}
\newtheorem{corollary}[theorem]{Corollary}
\newtheorem{algorithm}[theorem]{Algorithm}
\theorembodyfont{\rmfamily}
\newtheorem{definition}[theorem]{Definition}
\newtheorem{remark}[theorem]{Remark}
\newtheorem{example}[theorem]{Example}
\theorembodyfont{\rmfamily}

\newenvironment{proof}[1][\it Proof.]{\begin{trivlist}\item[\hskip \labelsep {\bfseries #1}]}{$\Box$\end{trivlist}}

\def\N{\ensuremath{\mathbb{N}}}
\def\Z{\ensuremath{\mathbb{Z}}}
\def\Q{\ensuremath{\mathbb{Q}}}
\def\R{\ensuremath{\mathbb{R}}}

\def\A{\ensuremath{\mathcal{A}}}
\def\F{\ensuremath{\mathcal{F}}}
\def\G{\ensuremath{\mathcal{G}}}

\def\C{\ensuremath{\mathcal{C}}}

\def\x{\ensuremath{{\bf{x}}}}

\def\init{\ensuremath{{\textup{in}}}}
\def\span{\ensuremath{{\textup{span}}}}
\def\flip{\ensuremath{{\textup{flip}}}}
\def\face{\ensuremath{{\textup{face}}}}

\def\vara{\ensuremath{x}}
\def\varb{\ensuremath{y}}
\def\varc{\ensuremath{z}}

\def\TAB{\hspace*{0.5cm}}

\begin{document}
\def\name{Gfan }

\title{Computing Gr\"obner fans}
\author{
Komei Fukuda\thanks{Research supported by the Swiss National Science Foundation Project  200021-105202,``Polytopes, Matroids and Polynomial Systems''.},
Anders N. Jensen\thanks{Research partially supported by the Faculty of Science, University of Aarhus, Danish Research Training Council (Forskeruddannelsesr\aa det, FUR) , Institute for Operations Research ETH, the Swiss National Science Foundation Project  200021-105202,
grants DMS 0222452 and DMS 0100141 of the U.S. National Science Foundation and the American Institute of Mathematics.}, and
Rekha R. Thomas\thanks{Research supported by grant DMS 0401047 of the U.S. National Science Foundation.}}
\maketitle

\begin{abstract}
  
  This paper presents algorithms for computing the Gr\"obner fan of an
  arbitrary polynomial ideal. The computation involves enumeration of
  all reduced Gr\"obner bases of the ideal. Our algorithms are based
  on a uniform definition of the Gr\"obner fan that applies to both
  homogeneous and non-homogeneous ideals and a proof that this object
  is a polyhedral complex. We show that the cells of a Gr\"obner fan
  can easily be oriented acyclically and with a unique sink, allowing
  their enumeration by the memory-less reverse search procedure. The
  significance of this follows from the fact that Gr\"obner fans are
  not always normal fans of polyhedra in which case reverse search
  applies automatically. Computational results using our
  implementation of these algorithms in the software package Gfan are
  included.

\end{abstract}
\normalsize

\section{Introduction}
The {\em Gr\"obner fan} of an ideal $I\subseteq k[x_1,\dots,x_n]$ was
defined by Mora and Robbiano in \cite{MoRo}. It is a fan of polyhedral
cones indexing {\em initial ideals} of $I$.
The full-dimensional cones are in bijection with the distinct monomial
initial ideals with respect to term orders or equivalently, the
reduced Gr\"obner bases of the ideal.  In this paper we will describe
algorithms for computing Gr\"obner fans of arbitrary polynomial ideals
allowing us to study their structure in detail. Our algorithms are
implemented in the software package \name \cite{gfan}.

The computation of the Gr\"obner fan of $I$ in terms of reduced
Gr\"obner bases yields \emph{a universal Gr\"obner basis} of $I$, a
set of polynomials which is a Gr\"obner basis of $I$ with respect to
every term order. The Gr\"obner fan also plays an important role in
{\em Gr\"obner basis conversion} \cite{collart} and the emerging field
of {\em tropical mathematics} as it contains the \emph{tropical
  variety} of $I$ as a subfan \cite{ctv}. Many of the well-known
theoretical applications of Gr\"obner bases rely on the existence of a
Gr\"obner basis or initial ideal for an ideal with prescribed
properties such as a particular complexity (a specified degree or
squarefree-ness) or homological properties (Cohen-Macaulay, Gorenstein
etc). The (even partial) computation of the Gr\"obner fan makes such
experimentations possible. No software package for the computation of
Gr\"obner fans is available at present.

In the literature a distinction is often made between the case of $I$
being homogeneous where the Gr\"obner fan is a complete fan in
$\R^n$, and the case of $I$ being non-homogeneous, where the
\emph{restricted} Gr\"obner fan with support equal to $\R_{\geq 0}^n$
is considered.  In this paper we avoid this distinction by giving a
new uniform definition of the Gr\"obner fan and we prove that this
actually defines a fan in the sense of a \emph{polyhedral complex} ---
a proof that was left out in \cite{MoRo} but was proven for the
special case of homogeneous ideals in \cite{sturmfels}. See
Section~\ref{section: groebner fan}.

In \cite[Algorithm 3.2 and 3.6]{sturmfels} and \cite[Section 6]{MoRo}
methods for computing the Gr\"obner fan of a polynomial ideal were
given. In this paper we shall study \cite[Algorithm 3.6]{sturmfels} in
detail. This algorithm traverses the maximal cones of the Gr\"obner
fan.
In the special case of \emph{toric} ideals the traversal
algorithm was already studied and implemented in
\cite{huber}.  The traversal is graph-like --- given a maximal cone
we need to be able to find its facets and we need to
be able to walk through a facet to the neighboring
maximal cone. Algorithms for performing these local
computations are discussed in Section \ref{section: computation}. These
amount to solving linear programming problems and using the local
basis change procedure due to \cite{collart}. We explain how to apply
these methods to our case.

The Gr\"obner fan of a homogeneous ideal $I$ is known to be the normal
fan of a polytope, the \emph{state polytope} of $I$ (\cite[Theorem
2.5]{sturmfels}). In this homogeneous case traversal of the maximal
cones in the Gr\"obner fan by walking through facets is equivalent to
traversal of the edge graph of the state polytope. In \cite{af-rse-96}
the memory-less \emph{reverse search} procedure for traversing the
edge graph of a polytope was given. This procedure easily applies to
Gr\"obner fans of homogeneous ideals. However, the question is what
happens if the ideal is not homogeneous. In \cite{jensen} a
\emph{non-regular} Gr\"obner fan was presented --- a fan that is not
the normal fan of any polyhedron. In light of this example it is not
clear that the reverse search technique applies to Gr\"obner fans in
general. In Section \ref{section: reverse search} we prove that all
Gr\"obner fans have what we shall call the \emph{reverse search
  property}, allowing them to be traversed efficiently.

Gr\"obner fans are often computed for ideals that possess a great deal
of symmetry.
In Section \ref{section: symmetry} we describe how to take advantage
of symmetry in the computations. The methods used here are similar to
those in Rambau's software package TOPCOM \cite{rambau} for
traversing the secondary fan of a point configuration up to symmetry.

In Section~\ref{section: complexity} we discuss the complexity of our
enumeration algorithm and in Section~\ref{sec experiments} we present
several examples of Gr\"obner fans computed using Gfan. This software
package uses the GNU multi-precision library \cite{gmp} for exact
arithmetics and Cddlib \cite{cdd} for solving linear programming
problems.

\section{The Gr\"obner fan of a polynomial ideal}
\label{section: groebner fan}

Let $R=k[x_1,\dots,x_n]$ be the polynomial ring in $n$ variables over
a field $k$ and let $I\subseteq R$ be an ideal. The \emph{Gr\"obner fan} and the \emph{restricted Gr\"obner fan} of $I$ are $n$-dimensional polyhedral fans defined in \cite{MoRo}. 
We recall the definition of a \emph{fan} in $\R^n$. A \emph{polyhedron} in $\R^n$ is a set of the form $\{x\in\R^n:Ax\leq b\}$ where $A$ is a matrix and $b$ is a vector. Bounded polyhedra are called \emph{polytopes}. If $b=0$ the set is a \emph{polyhedral cone}.
The \emph{dimension} of a polyhedron is the dimension of the smallest affine subspace containing it. A \emph{face} of a polyhedron $P$ is either the empty set or a non-empty subset of $P$ which is the set of maximizers of a linear form over $P$. We use the following notation for the face maximizing a form $\omega\in\R^n$:
$$ \face_\omega(P)=\{p\in P:\langle \omega,p\rangle=\textup{max}_{q\in P}\langle \omega,q\rangle\}.$$
A face of $P$ is called a \emph{facet} if its dimension is one smaller than the dimension of $P$.  
\begin{definition}
A collection $\C$ of polyhedra in $\R^n$ is a \emph{polyhedral complex} if:
\begin{enumerate}
\item all non-empty faces of a polyhedron $P\in \C$ are in $\C$, and
\item the intersection of any two polyhedra $A,B\in \C$ is a face of $A$ and a face of $B$. 
\end{enumerate}
The \emph{support} of $\C$ is the union of its members. A polyhedral complex is a \emph{fan} if it only consists of cones. A fan is \emph{pure} if all its maximal cones have the same dimension.
\end{definition}
A simple way to construct a fan is to take the \emph{normal fan} of a
polyhedron.   
\begin{definition}
Let $P\subseteq\R^n$ be a polyhedron.
For a face $F$ of $P$ we define its \emph{normal cone}
$$N_P(F):=\overline{\{\omega\in\R^n:\face_\omega(P)=F\}}$$
with the closure being taken in the usual topology. The \emph{normal fan} of $P$ is the fan consisting of the normal cones $N_P(F)$ as $F$ runs through all non-empty faces of $P$.
\end{definition}
If the support of a fan is $\R^n$, the fan is said to be
\emph{complete}. It is clear that the normal fan of a polytope is
complete.  Not all fans arise as the normal fan of a polyhedron
\cite[page 25]{Ful}.

For $\alpha\in\N^n$ we use the notation $\x^\alpha:=x_1^{\alpha_1}\dots x_n^{\alpha_n}$ for a monomial in $R$. By
a \emph{term order} on $R$ we mean a total ordering on all monomials in $R$
such that:
\begin{enumerate}
\item For all $\alpha\in\N^n\backslash \{0\}:1<\x^\alpha$ and
\item for $\alpha,\beta,\gamma\in\N^n: \x^\alpha < \x^\beta \Rightarrow \x^\alpha \x^\gamma < \x^\beta \x^\gamma$.
\end{enumerate}
By a {\em term} we mean a monomial together with its coefficient.
Term orders are used for ordering terms, ignoring the coefficients.
For a vector $\omega\in\R_{\geq 0}^n$ and a term order $\prec$ we
define the new term order $\prec_\omega$ as follows:
$$\x^\alpha\prec_\omega \x^\beta ~\Longleftrightarrow~
\langle\omega,\alpha\rangle < \langle\omega,\beta\rangle ~\vee~
(\langle\omega,\alpha\rangle =\langle\omega,\beta\rangle ~\wedge~
\x^\alpha\prec \x^\beta).$$
Let $\prec$ be a term order. For a
non-zero polynomial $f\in R$ we define its \emph{initial term},
$\init_\prec(f)$, to be the unique maximal term of $f$ with respect to
$\prec$.  In the same way for $\omega\in \R^n$ we define the
\emph{initial form}, $\init_\omega(f)$, to be the sum of all terms of
$f$ whose exponents maximize $\langle \omega,\cdot\rangle$. The
polynomial $f$ is $\omega$\emph{-homogeneous} if $\init_\omega(f)=f$.
The $\omega$\emph{-degree} of a term $c \x^\alpha$ is $\langle
\omega,\alpha\rangle$ and the $\omega$\emph{-degree} of a non-zero
polynomial $f$ is the common $\omega$-degree of the terms of
$\init_\omega(f)$. The \emph{initial ideals} of an ideal $I$ with
respect to $\prec$ and $\omega$ are defined as
$$ \init_\prec(I)=\langle \init_\prec(f):f\in I\backslash\{0\}\rangle \,\,\mbox{and}\,\, \init_\omega(I)=\langle \init_\omega(f):f\in I\rangle.$$
Note that $\init_\prec(I)$ is a monomial ideal while $\init_\omega(I)$ might not be. A monomial in $R\backslash \init_\prec(I)$ (with coefficient $1$) is called a \emph{standard monomial} of $\init_\prec(I)$.

Although initial ideals are defined with respect to not necessarily positive vectors, Gr\"obner bases are only defined with respect to true term orders:
\begin{definition}
\label{def: groebner basis}
Let $I\subseteq R$ be an ideal and $\prec$ a term order on $R$. A generating set $\G=\{g_1,\dots,g_m\}$ for $I$ is called a \emph{Gr\"obner basis} for $I$ with respect to $\prec$ if
$$\init_\prec(I)=\langle \init_\prec(g_1),\dots,\init_\prec(g_m)\rangle.$$
The Gr\"obner basis $\G$ is \emph{minimal} if $\{\init_\prec(g_1),\dots,\init_\prec(g_m)\}$ generates $\init_\prec(I)$ minimally. A minimal Gr\"obner basis is \emph{reduced} if the initial term of every $g\in\G$ has coefficient $1$ and all other monomials in $g$ are standard monomials of $\init_\prec(I)$.
\end{definition}
We use the term \emph{marked Gr\"obner basis} for a Gr\"obner basis where the initial terms have been distinguished from the non-initial ones (they have been marked). For example, $\{\underline{x^2}+xy+y^2\}$ and $\{x^2+xy+\underline{y^2}\}$ are marked Gr\"obner bases for the ideal $\langle x^2+xy+y^2\rangle$ while $\{x^2+\underline{xy}+y^2\}$ is not since $xy$ is not the initial term of $x^2+xy+y^2$ with respect to any term order.

For a term order $\prec$ and an ideal $I$, Buchberger's algorithm
guarantees the existence of a unique marked reduced Gr\"obner basis.
We denote it by $\G_\prec(I)$.
For two term orders $\prec$ and $\prec'$, if $\init_\prec(I)=\init_{\prec'}(I)$ then $\G_\prec(I)=\G_{\prec'}(I)$. Conversely, given a marked Gr\"obner basis $\G_\prec(I)$, $\init_\prec(I)$ can be easily read off. 

Given an ideal $I$, a natural equivalence relation on $\R^n$ is
induced by taking initial ideals:
\begin{equation}
\label{eq:eq}
u\sim v ~\Longleftrightarrow~ \init_u(I)=\init_v(I).
\end{equation}
We introduce the following notation for the closures of the
equivalence classes: 
$$C_\prec(I)=\overline{\{u\in\R^n : \init_u(I)=\init_\prec(I)\}}\mbox{~~and}$$
$$C_v(I)=\overline{\{u\in\R^n : \init_u(I)=\init_v(I)\}}.~~~~~$$
\begin{remark}
\label{rem: facts}
It is well known that for a fixed ideal $I$ there are only finitely
many sets $C_\prec(I)$ and they cover $\R_{\geq 0}^n$, see
\cite{MoRo}. Secondly, every initial ideal $\init_\prec(I)$ is of the
form $\init_\omega(I)$ for some $\omega\in\R_{>0}^n$, see
\cite[Proposition 1.11]{sturmfels}. Consequently, every $C_\prec(I)$
is of the form $C_\omega(I)$.
\end{remark}
A third observation is that the equivalence classes are not convex in
general since we allow the vectors to be anywhere in $\R^n$:
\begin{example}
\label{first example}
Let $I=\langle \vara-1, \varb-1\rangle$. The ideal $I$ has five initial
ideals: $\langle \vara-1, \varb-1\rangle$, $\langle \vara, \varb\rangle$,
$\langle \vara, \varb-1\rangle$, $\langle \vara-1, \varb\rangle$ and $\langle
1\rangle$. In particular, for $u=(-1,3)$ and $v=(3,-1)$ we have
$\init_u(I)=\init_v(I)=\langle 1\rangle$ but $\init_{{1\over
2}(u+v)}(I)=\langle \vara, \varb \rangle$.
\end{example}
\begin{proposition}
\label{l4}
\label{prop: groebner cone}
Let $\prec$ be a term order and $v\in C_\prec(I)$. For
$u\in\R^n$
$$\init_u(I)=\init_v(I) ~\Longleftrightarrow ~\forall
g\in\G_\prec(I),\,\init_u(g)=\init_v(g).$$  
\end{proposition}
This proposition is a little more general than Proposition 2.3 in
\cite{sturmfels} as it allows the vectors $u$ and $v$ to have negative
components. A proof is given in the next section. For fixed $\prec$
and $v$ as in Proposition~\ref{prop: groebner cone}, we get that
$C_v(I)$, the closure of the equivalence class of $v$, is a polyhedral
cone since each $g\in\G_\prec(I)$ introduces the equation
$\init_u(g)=\init_v(g)$ which is equivalent to having $u$ satisfy a
set of linear equations and strict linear inequalities, see Example
\ref{ex:gfanbig}. The closure is obtained by making the strict
inequalities non-strict.  Under the assumptions of Proposition
\ref{prop: groebner cone} we may write this in the following way:
\begin{equation}
\label{eq:eq2}
u\in C_v(I) ~\Longleftrightarrow ~\forall
g\in\G_\prec(I),\,\init_v(\init_u(g))=\init_v(g). 
\end{equation}
As we saw in Example~\ref{first example}, not all equivalence classes
are convex. However, for an arbitrary $v$, $C_v(I)$ is a convex
polyhedral cone if it contains a strictly positive vector. In this
case, there must exist a vector $p\in\R_{>0}^n$ in the interior of
$C_v(I)$. Then $\init_p(I)=\init_v(I)$ and, by Lemma \ref{l6.1}, $p\in
C_{\prec_p}(I)$ for any $\prec$. Hence the equivalence class of $v$ is
of the form required in Proposition \ref{prop: groebner cone}.

\begin{figure}
\label{fig:gfanbig}
\begin{center}
\epsfig{file=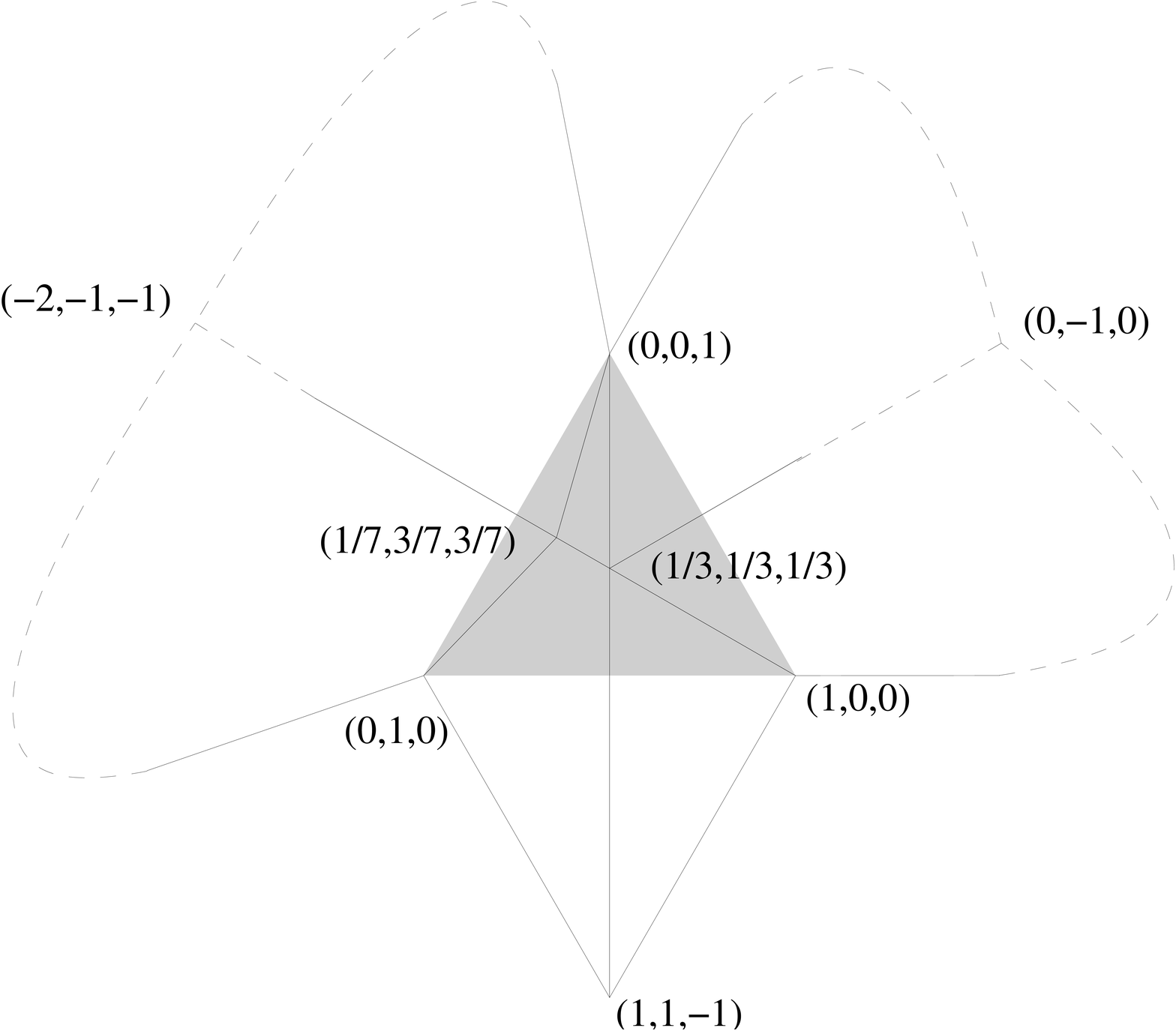,width=7cm}
\end{center}
\caption{The Gr\"obner fan of the ideal in Example \ref{ex:gfanbig} has 7 three-dimensional, 14 two-dimensional and 8 one-dimensional cones. The intersections of the two-dimensional cones with the hyperplane $\vara+\varb+\varc=1$ are drawn as lines. The dotted part of the figure shows the combinatorial structure outside the hyperplane. The gray triangle indicates the positive orthant.}
\end{figure}
\begin{example}
\label{ex:gfanbig}
Let $I=\langle
\vara+\varb+\varc,\vara^3\varc+\vara+\varb^2\rangle\subseteq
\Q[\vara,\varb,\varc]$ and let $\prec$ be the lexicographic term order
with $\vara\prec \varb\prec \varc$.  Then
$\G_\prec(I)=\{\underline{\varb^2}+\vara-\vara^3\varb-\vara^4,\underline{\varc}+\varb+\vara\}$.
If $v=(1,4,5)$ then $\init_v(I)=\init_\prec(I)=\langle
\varb^2,\varc\rangle$ and $C_v(I)=C_\prec(I)$. By Proposition
\ref{l4}, $\init_u(I)=\init_v(I)$ if and only if the following two
equations are satisfied:
$$\init_u(\varc+\varb+\vara)=\varc~(\Leftrightarrow~ u_\varc>\textup{max}\{u_\vara,u_\varb\}),\textup{ and}$$
$$\init_u(\varb^2+\vara-\vara^3\varb-\vara^4)=\varb^2~(\Leftrightarrow~ 2u_\varb>\textup{max}\{u_\vara,3u_\vara+u_\varb,4u_\vara\}).$$
Introducing non-strict inequalities we obtain a description of $C_\prec(I)$. This cone is simplicial and has the cones $C_{(0,0,1)}(I)$, $C_{(1,3,3)}(I)$ and $C_{(-2,-1,-1)}(I)$ as extreme rays and $C_{(1,3,4)}(I)$, $C_{(-2,-1,0)}(I)$ and $C_{(-1,2,2)}(I)$ as facets. Since $(-2,-1,0)$ is in $C_\prec(I)$ a description of vectors $u$ in $C_{(-2,-1,0)}(I)$ is given by:
$$\init_{(-2,-1,0)}(\init_u(\varc+\varb+\vara))=\varc ~(\Leftrightarrow~ u_\varc\geq\textup{max}\{u_\vara,u_\varb\}),\textup{ and}$$
$$\init_{(-2,-1,0)}(\init_u(\varb^2+\vara-\vara^3\varb-\vara^4))=\varb^2+\vara ~(\Leftrightarrow~ 2u_\varb=u_\vara\geq\textup{max}\{3u_\vara+u_\varb,4u_\vara\}).$$
\end{example}

\begin{figure}
\label{fig:gfan}
\begin{center}
\epsfig{file=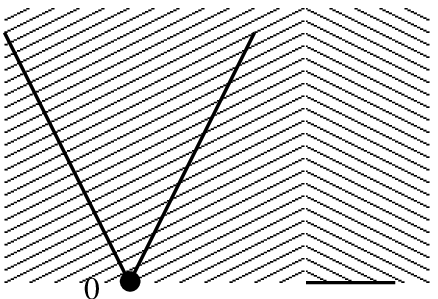,height=3.4cm}
\epsfig{file=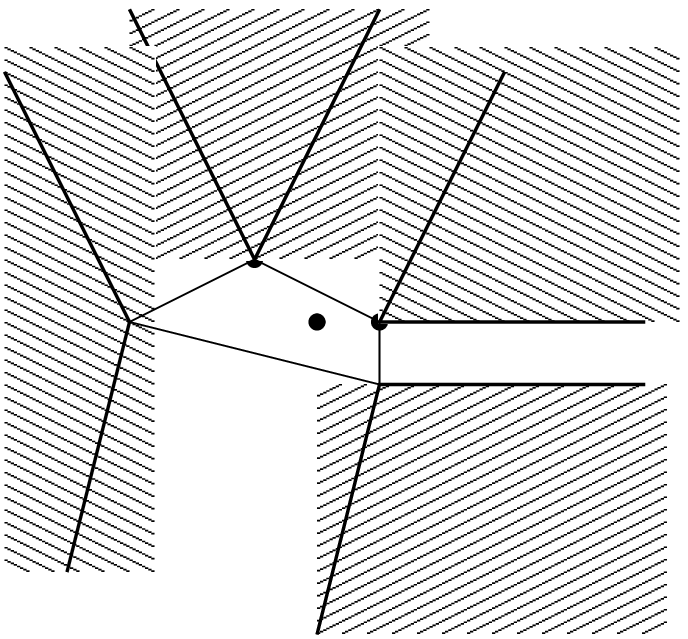,height=3.4cm}~~
\epsfig{file=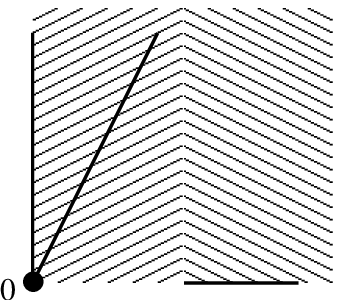,height=3.4cm}
\end{center}
\caption{The Gr\"obner fan of the ideal in Example \ref{ex:gfan} is
  shown on the left. The restricted Gr\"obner fan is on the right. In
  the middle the Newton polytope of the generator is drawn with the
  shape of its normal fan indicated.} 
\end{figure}

\begin{definition}
\label{definition gfan}
The \emph{Gr\"obner fan} of an ideal $I\subseteq R$ is the set of the
closures of all equivalence classes intersecting the positive orthant
together with their proper faces.
\end{definition}
This is a variation of the definitions appearing in the
literature. The advantage of this variant is that it gives
well-defined and \emph{nice} fans in the homogeneous and
non-homogeneous case simultaneously. By \emph{nice} we mean that all cones in this fan are closures of equivalence classes.
It is not clear a priori that the Gr\"obner fan is a polyhedral complex. A proof is given in the next section (Theorem \ref{thm: gfan is a fan}).  The support of the Gr\"obner fan of $I$ is called the \emph{Gr\"obner region} of $I$.
Recall that the common refinement of two fans $\F_1$ and $\F_2$ in $\R^n$ is defined as
$$ \F_1\wedge\F_2=\{C_1\cap C_2\}_{(C_1,C_2)\in\F_1\times\F_2}.$$
The common refinement of two fans is a fan. We define the \emph{restricted} Gr\"obner fan of an ideal to be the common refinement of the Gr\"obner fan and the faces of the non-negative orthant. The support of the restricted Gr\"obner fan is $\R_{\geq 0}^n$.
The \emph{Newton polytope} of a polynomial is the convex hull of its exponent vectors.
\begin{example}
\label{ex:gfan}
The Gr\"obner fan of the principal ideal $\langle \vara^4+\vara^4\varb-\vara^3\varb+\vara^2\varb^2+\varb\rangle$ consists of one $0$-dimensional cone, three $1$-dimensional cones and two $2$-dimensional cones, see Figure \ref{fig:gfan}. The same is true for the restricted Gr\"obner fan. Notice, however, that in the restricted Gr\"obner fan one of the $1$-dimensional cones and one of the $2$-dimensional cones are not equivalence classes of the equivalence relation (\ref{eq:eq}).
\end{example}

\subsection{Proof: The Gr\"obner fan is a fan}
In this section we prove that the Gr\"obner fan is a fan i.e., that it
is a polyhedral complex consisting of cones.  Recall, in general the
Gr\"obner fan is not complete and its support is larger than $\R_{\geq
  0}^n$. In \cite{MoRo} there is no proof that the Gr\"obner fan is a
fan in the sense of a polyhedral complex. A proof that the Gr\"obner
fan is a polyhedral complex under the assumption that the ideal is
homogeneous is given in \cite{sturmfels}.  We present a complete proof
for the general case. Many of the results we need in the proof are
generalizations of known results needed in the proof that the
Gr\"obner fan of a homogeneous ideal is a polyhedral complex
\cite{sturmfels}. However, we do not rely on these references for the 
sake of being self-contained. 

We fix the ideal $I\subseteq R$ in the following theorems. The most
important step is the proof of Proposition \ref{l4} which tells us
that the closure of an equivalence class is a polyhedral cone. Then we
prove that the relative interior of any face in the Gr\"obner fan is
an equivalence class (Proposition \ref{l6}) and, finally, that the
intersection of two cones in the fan is a face of both
(Proposition \ref{l9}). 

To prove Proposition \ref{l4} we start by proving a similar statement
for the equivalence classes arising from initial ideals with respect
to term orders.

\begin{lemma}
\label{l1}
Let $\prec$ be a term order. For $v\in\R^n$,
$$\init_v(I)=\init_\prec(I)~\Longleftrightarrow~ \forall g\in\G_\prec(I),
\,\init_v(g)=\init_\prec(g).$$
\end{lemma}
\begin{proof}
  $\Rightarrow$: Let $g\in\G_\prec(I)$. Since $\G_\prec(I)$ is
  reduced, only one term from $g$, $\init_\prec(g)$, can be in
  $\init_\prec(I)=\init_v(I)$. The initial ideal $\init_v(I)$ is a
  monomial ideal, implying that all terms of an element in the ideal
  must be in the ideal too. Hence, the initial form $\init_v(g)\in
  \init_v(I)$ has to be equal to $\init_\prec(g)$.

\noindent
$\Leftarrow$: We must show that $\init_v(I)=\init_\prec(I)$ where
$\init_\prec(I) = \langle \init_\prec(g)\rangle_{g\in\G_\prec(I)}$.
The ``$\supseteq$'' inclusion is clear since
$\init_\prec(g)=\init_v(g) \in \init_v(I)$ for all $g\in\G_\prec(I)$.

To prove the ``$\subseteq$'' inclusion, since $\init_v(I) = \langle
\init_v(f), \, f\in I \rangle$, it suffices to show that $\init_v(f)
\in \init_\prec(I)$ for all $f \in I$. Pick $f\in I$ and reduce it to
zero using the division algorithm (e.g.  \cite[Chapter 2]{Cox96}) with
$\G_\prec(I)$ and $\prec$. We may write
\begin{equation}
\label{e1}
f= m_1g_{i_1}+\dots +m_rg_{i_r}
\end{equation}
where $m_j$ is a monomial and $g_{i_j}$ is an element from
$\G_\prec(I)$. The division algorithm guarantees that
$\init_\prec(f)\ge m_j \init_\prec(g_{i_j})$ with respect to $\prec$
since monomials are substituted with monomials less than the original
ones with respect to $\prec$ in the division process. Exactly the same
thing is true for $v$-degrees since $v$ and $\prec$ agree on
$\G_\prec(I)$. Thereby, any monomial on the right hand side in
(\ref{e1}) has $v$-degree less than or equal to the $v$-degree of the
left hand side.  Consequently,
$$\init_v(f)=\sum_{j\in J}m_j \init_v(g_{i_j})$$
with $j$ running
through a subset such that $m_j \init_v(g_{i_j})$ has the same
$v$-degree as $\init_v(f)$. Since $\init_v(g)\in \init_\prec(I)$, the
initial form $\init_v(f)\in \init_\prec(I)$.
\end{proof}

By Lemma~\ref{l1} the equivalence class of $\init_\prec(I)$ is open.
Since $\init_\prec(I)$ is of the form $\init_v(I)$ for some $v$ (see
Remark \ref{rem: facts}), the equivalence class of $\init_\prec(I)$ is
also non-empty and hence full-dimensional. Thus we have proved that
the equivalence class of a term order is a full dimensional open
polyhedral cone.

\begin{corollary}
\label{l2.1}
Let $\prec$ be a term order and $v\in\R^n$. Then
$$v\in{C_\prec(I)} \Leftrightarrow \forall g\in\G_\prec(I):
\init_\prec(\init_v(g))=\init_\prec(g).$$ 
\end{corollary}
\begin{proof}
  Lemma~\ref{l1} tells us that $v$ lies in the interior of
  $C_\prec(I)$ if and only if $\init_v(g) = \init_\prec(g)$ for all $g
  \in \G_\prec(I)$. Relaxing the resulting strict inequalities to
  non-strict inequalities we get a description of $C_\prec(I)$. This
  relaxation is exactly the one given by
  $\init_\prec(\init_v(g))=\init_\prec(g)$ for all $g$ in
  $\G_\prec(I)$.
\end{proof}

\begin{lemma}
\label{l2.2}
A polynomial $f\in \init_v(I)$ can be written in the form $f=\sum_i
\init_v(c_i)$ where $c_i\in I$ and all summands in the sum have
different $v$-degrees.
\end{lemma}
\begin{proof}
  The initial ideal $\init_v(I)$ is generated by $v$-homogeneous
  polynomials, implying that all $v$-homogeneous components of $f$ are
  in $\init_v(I)$. Let $h$ be a maximal $v$-homogeneous component of
  $f$.  We need to show that $h$ is the initial form of an element in
  $I$ with respect to $v$. We may write $h$ as $\init_v(a_1)+\dots
  +\init_v(a_s)$ for some polynomials $a_1,\dots,a_s$ in $I$. Since
  $h$ is $v$-homogeneous we can rewrite $h$ as the sum $\sum_{j\in
    J}\init_v(a_j)$ of forms having the same $v$-degree as $h$. We
  pull out the initial form and get $h=\init_v(\sum_{j\in J}a_j)$.
\end{proof}

\begin{lemma}
\label{l2}
Let $\prec$ be a term order. If $v\in {C_\prec(I)}$ then
$\init_\prec(\init_v(I))=\init_\prec(I)$.
\end{lemma}
\begin{proof}
  Let $g\in\G_\prec(I)$. Since $v\in {C_\prec(I)}$, by Corollary
  \ref{l2.1}, $\init_\prec(g)=\init_\prec(\init_v(g))$ and hence
  $\init_\prec(I) = \langle \init_\prec(g)\rangle_{g\in\G_\prec(I)}
  \subseteq \init_\prec(\init_v(I))$.
  
  We now prove that $\init_\prec(\init_v(I)) \subseteq
  \init_\prec(I)$. Notice that $\init_\prec(\init_v(I))$ is generated
  by initial terms of elements $f \in \init_v(I)\backslash\{0\}$ with respect to
  $\prec$.  Suppose $f\in \init_v(I)\backslash\{0\}$. It suffices to show that
  $\init_\prec(f) \in \init_\prec(I)$. Using Lemma \ref{l2.2} we may
  write $f=\sum_{i=1}^s \init_v(c_i)$ where $c_1,\dots,c_s\in I$ and
  $\init_v(c_1),\dots,\init_v(c_s)$ are $v$-homogeneous each with
  distinct degree, so that no cancellations occur. Consequently
  $\init_\prec(f)$ equals $\init_\prec(\init_v(c_j))$ for some $j$.
  We wish to prove that $\init_\prec(\init_v(c_j)) \in
  \init_\prec(I)$.  We use the division algorithm with $\G_\prec(I)$
  and $\prec$ to rewrite $c_j$
  $$c_j=m_1g_{i_1}+\dots+m_rg_{i_r}$$
  where $m_1,\dots,m_r$ are
  monomials and $g_{i_1},\dots, g_{i_r}$ belong to $\G_\prec(I)$.  Let
  $M$ be the $v$-degree of $c_j$. In the division algorithm we
  sequentially reduce $c_j$ to zero. In each step, the $v$-degree of
  $c_j$ will decrease or stay the same since we subtract the product
  of a monomial and an element from $\G_\prec(I)$ where the $v$-degree
  of the product already appeared in $c_j$ by Corollary \ref{l2.1}.
  Equivalently, the product of the monomial and the element from
  $\G_\prec(I)$ are ``added'' to the right hand side of the equation.
  We are done when $c_j=0$ and or equivalently, the original $c_j$ is
  written as the above sum with every term having $v$-degree less or
  equal to $M$.  Consequently, we have
  $$\init_v(c_j)=\sum_{j'\in J'}\init_v(m_{j'}g_{i_{j'}})$$
  for a
  suitable $J'$.  The division algorithm guarantees that the exponent
  vectors of $\init_\prec(m_1g_{i_1}),\dots,\init_\prec(m_rg_{i_r})$
  are distinct. Since $v \in C_\prec(I)$, they equal
  $\init_\prec(\init_v(m_1g_{i_1})),\dots,\init_\prec(\init_v(m_rg_{i_r}))$.
  The maximal one of these with respect to $\prec$ cannot cancel in
  the sum. Hence
  $\init_\prec(\init_v(c_j))=\init_\prec(m_{j'}g_{i_{j'}})$ for some
  $j'$ which implies that  $\init_\prec(\init_v(c_j)) \in
  \init_\prec(I)$ as needed.
\end{proof}

An easy corollary is a method for computing Gr\"obner bases for initial ideals.

\begin{corollary}
\label{l3}
Let $\prec$ be a term order. If $v\in {C_\prec(I)}$ then
$$\G_\prec(\init_v(I))=\{\init_v(g)\}_{g\in\G_\prec(I)}.$$
\end{corollary}

\begin{proof}
  By Corollary \ref{l2.1}, $\langle
  \init_\prec(\init_v(g))\rangle_{g\in\G_\prec(I)}=\langle
  \init_\prec(g)\rangle_{g\in\G_\prec(I)}=\init_\prec(I)$. By Lemma
  \ref{l2}, $\init_\prec(I)$ equals $\init_\prec(\init_v(I))$. Thus
  $\init_\prec(\init_v(I)) = \langle
  \init_\prec(\init_v(g))\rangle_{g\in\G_\prec(I)}$. This proves that
  $\{\init_v(g)\}_{g\in\G_\prec(I)}$ is a Gr\"obner basis of
  $\init_v(I)$ with respect to $\prec$. It is reduced since
  $\G_\prec(I)$ is minimal and reduced.
\end{proof}

\noindent
We are now able to give a proof for Proposition \ref{prop: groebner
  cone} which claimed that given $v \in C_\prec(I)$ and $u\in\R^n$,
$\init_u(I)=\init_v(I) ~\Longleftrightarrow ~\forall
g\in\G_\prec(I),\,\init_u(g)=\init_v(g).$
  
\begin{proof}
  $\Leftarrow:$ Since $\init_u(g)=\init_v(g)$ for all
  $g\in\G_\prec(I)$, we get that $\init_\prec(\init_u(g)) =
  \init_\prec(\init_v(g))$ for all $g\in\G_\prec(I)$. Since $v\in
  {C_\prec(I)}$, by Corollary \ref{l2.1}, $\init_\prec(g) =
  \init_\prec(\init_v(g))$ for all $g \in \G_\prec(I)$ and hence
  $\init_\prec(g) = \init_\prec(\init_u(g))$ for all $g \in
  \G_\prec(I)$ and $u\in {C_\prec(I)}$ by Corollary \ref{l2.1}. The Gr\"obner basis
  $\G_\prec(\init_u(I))$ is then $\{ \init_u(g)\}_{g\in\G_\prec(I)}$
  by Corollary \ref{l3}. We get the same Gr\"obner basis for
  $\init_v(I)$. Hence, $\init_u(I)=\init_v(I)$.

\noindent
$\Rightarrow$: Let $g\in\G_\prec(I)$. We need to show that $\init_u(g)
= \init_v(g)$. Since the basis is reduced, only one term of $g$, namely
$\init_\prec(g)$, is in $\init_\prec(I)$.  We start by proving that
the term $\init_\prec(g)$ is a term in $\init_v(g)$ and a term in
$\init_u(g)$. For $\init_v(g)$ we apply Corollary \ref{l2.1} which
says $\init_\prec(g)=\init_\prec(\init_v(g))$. For $\init_u(g)$ we
apply Lemma \ref{l2} and get $\init_\prec(\init_u(g))\in
\init_\prec(\init_u(I))=\init_\prec(\init_v(I))=\init_\prec(I)$. Only
one term of $g$ is in $\init_\prec(I)$, so
$\init_\prec(\init_u(g))=\init_\prec(g)$.  If the difference
$\init_u(g)-\init_v(g)$, belonging to $\init_u(I)=\init_v(I)$, is
non-zero we immediately reach a contradiction since the difference
contains no terms from $\init_\prec(I)=\init_\prec(\init_v(I))$.
\end{proof}

We have now proved that every equivalence class of a vector $v$ in a
$C_\prec(I)$ is a relatively open convex polyhedral cone. By the argument following Proposition \ref{l4} in the previous section all sets in the Gr\"obner fan are in fact cones.
We now argue that the relative interior of every cone
in the Gr\"obner fan is an equivalence class. 

\begin{lemma}
\label{l6.1}
Let $\prec$ be a term order. If $v\in\R_{\geq 0}^n$ then
$v\in{C_{\prec_v}(I)}$.
\end{lemma}
\begin{proof}
  This follows from Corollary \ref{l2.1} since
  $\init_{\prec_v}(\init_v(g))=\init_{\prec_v}(g)$ for all $g \in
  \G_{\prec_v}(I)$.
\end{proof}

\begin{proposition}
\label{l6}
The relative interior of a cone in the Gr\"obner fan is an equivalence
class (with respect to $u\sim {u'} \Leftrightarrow
\init_u(I)=\init_{u'}(I)$).
\end{proposition}
\begin{proof}
  By definition every cone in the fan is the face of the closure of an
  equivalence class for a positive vector $v\in\R_{>0}^n$. Let
  $\prec'$ be an arbitrary term order and define $\prec$ as
  $\prec'_v$. According to Lemma \ref{l6.1} the vector $v$ belongs to
  ${C_{\prec}(I)}$. Notice that by (\ref{eq:eq2}), $C_v(I)\subseteq
  C_\prec(I)$ since for all $u\in C_v(I)$ and $g\in\G_\prec(I)$, the
  condition $\init_\prec(\init_u(g))=\init_\prec(\init_v(\init_u(g)))=
  \init_\prec(\init_v(g))=\init_\prec(g)$ of Corollary \ref{l2.1} is
  satisfied.
  By (\ref{eq:eq2}) the closed set ${C_v(I)}$ is cut out by some
  equations and non-strict inequalities. The relative interior of any
  face of ${C_v(I)}$ can be formed from this inequality system by
  changing a subset of the inequalities to strict inequalities and the
  remaining ones to equations. So let $u$ be a vector in the relative
  interior of some face of ${C_v(I)}$. The vector $u$ is in
  $C_v(I)\subseteq{C_{\prec}(I)}$. We may use Proposition \ref{l4} to
  conclude that a vector $u'\in\R^n$ is equivalent to $u$ if and only
  if it satisfies the inequality system mentioned above --- that is,
  if and only if it is in the relative interior of the face.
\end{proof}

It remains to be shown that the intersection of two cones in the
Gr\"obner fan is a face of both cones (Proposition \ref{l9}). We need
a few observations.

\begin{corollary}
\label{l8}
Let $C$ be a cone in the Gr\"obner fan. If $v\in C$ then for
$u\in\R^n$, 
$$\init_u(I)=\init_v(I) ~\Rightarrow~ u\in C.$$
\end{corollary}

\begin{proof}
  The vector $v$ is in the relative interior of some face of ${C}$.
  This face is also in the Gr\"obner fan. By Proposition \ref{l6} $u$
  is in the relative interior of the same face and, consequently, also
  in ${C}$.
\end{proof}

By Remark \ref{rem: facts} there are only finitely many initial ideals
given by term orders and, consequently, only finitely many reduced
Gr\"obner bases of $I$. It follows that there can only be finitely
many equivalence classes of the type described in Proposition
\ref{prop: groebner cone} and Proposition \ref{l6}.

\begin{proposition}
\label{l9}
Let ${C_1}$ and ${C_2}$ be two cones in the Gr\"obner fan of $I$. Then
the intersection ${C_1}\cap{C_2}$ is a face of $C_1$.
\end{proposition}

\begin{proof}
  The intersection ${C_1}\cap{C_2}$ is a cone. By Corollary \ref{l8},
  $C_1$ and $C_2$ are unions of equivalence classes. Further, if $v
  \in C_1 \cap C_2$, then again by Corollary \ref{l8}, the entire
  equivalence class of $v$ is both in $C_1$ and in $C_2$ and hence in
  $C_1 \cap C_2$. Hence ${C_1}\cap{C_2}$ is a union of equivalence
  classes.
  
  Let $u$ be a vector in such an equivalence class $E$ contained in
  $C_1\cap C_2$. Then $u$ is in the relative interior of one of the
  faces of $C_1$ which is a cone in the Gr\"obner fan. By Proposition
  \ref{l6} the set of vectors in the relative interior of this face is
  exactly $E$. Hence every such equivalence class is the relative
  interior of a face of $C_1$ and its closure is the face.
  
  Look at the $\R$-span of each equivalence class contained in
  $C_1\cap C_2$. These spans must be different for every face of
  $C_1$. We claim that there can be only one maximal dimensional
  cone/span. If there were two cones then their convex hull would be
  in $C_1\cap C_2$ and have dimension at least one higher and thus
  cannot be covered by the finitely many lower dimensional equivalence
  classes --- a contradiction.

Let $E$ be the maximal dimensional equivalence class contained in
$C_1\cap C_2$. We will argue that $\overline{E}=C_1\cap C_2$. The
inclusion $\overline{E}\subseteq C_1\cap C_2$ is already clear since
$C_1\cap C_2$ is closed. To prove the other inclusion suppose
$\omega\in C_1\cap C_2\backslash \overline{E}$. Then
$\textup{conv}(\overline{E},\omega)\backslash\overline{E}$ is
contained in $C_1\cap C_2$ and has dimension at least the dimension of
$E$. This is a contradiction since
$\textup{conv}(\overline{E},\omega)\backslash\overline{E}$ cannot be
covered by finitely many lower dimensional equivalence classes. This
completes the proof.
\end{proof}

\begin{theorem}
\label{thm: gfan is a fan}
\label{l10}
The Gr\"obner fan is a polyhedral complex of cones and hence a fan.
\end{theorem}
\begin{proof}
  We already argued using Proposition \ref{l4} and Lemma \ref{l6.1}
  that the Gr\"obner fan consists of polyhedral cones.  The first
  condition for being a polyhedral complex is satisfied by definition.
  The second condition is Proposition \ref{l9}.
\end{proof}

\section{Reverse search property}
\label{section: reverse search}

By the \emph{graph} of a pure full-dimensional fan we mean the set of maximal cones with two cones being connected if they share a common facet. In this section we will prove that the \emph{reverse search technique} \cite{af-rse-96} can be used for traversing the graph of a Gr\"obner fan.
This follows from the main theorem, Theorem \ref{thm:gfan has reverse search property}, which says that the graph of a Gr\"obner fan can be oriented easily without cycles and with a unique sink.
 In Definition \ref{def:reverse search property} we define what we mean by this.

We start by explaining how a graph with this special kind of orientation can be traversed by reverse search. The idea is to define a spanning tree of the graph which can be easily traversed. The following is a simple proposition which we shall not prove.

\begin{proposition}
\label{proposition: reverse search graph}
Let $G=(V,E)$ be an oriented graph without cycles and with a unique
sink $s$. If for every vertex $v\in V\backslash\{s\}$ some outgoing
\emph{search edge} $e_v=(v,\cdot)$ is chosen then the set of chosen edges
is a spanning tree for $G$.
\end{proposition}
The spanning tree in Proposition \ref{proposition: reverse search graph} is referred to as the \emph{search tree}. The proposition implies that the graph is connected.

Notice that we can find the sink by starting at any vertex and walking
along a unique path of search edges until we get stuck, in which case
we are at the sink. Consequently, the sink is the root of the oriented
spanning tree.  A corollary to the proposition is the \emph{reverse
  search algorithm} for traversing $G$:
\begin{algorithm}
\label{alg: reverse search}  
Let $G=(V,E)$ be the oriented graph of Proposition \ref{proposition: reverse
search graph} and suppose the choice of a search edge $e_v$ for each
vertex $v\not=s$ has been made. Calling the following recursive
procedure with $v=s$ will output all vertices
in $G$.\\
{\bf Output\_subtree}(v)\\
{\bf Input:} A vertex $v$ in the graph $G$.\\
{\bf Output:}
The set of vertices in the subtree with root $v$.\\
$\{$\\
\TAB Output
$v$;\\
\TAB Compute the edges of form $(\cdot,v)\in E$;\\
\TAB For every oriented edge $(u,v)\in E$\\
\TAB \TAB If
($e_{u}=(u,v)$) {\bf Output\_subtree}(u);\\
$\}$
\end{algorithm}
This algorithm does not have to store a set of ``active'' vertices as
is usually needed in depth- and breadth-first traversals. It is
even possible to formulate the algorithm completely without recursion
avoiding the need for a recursion stack. In that sense the algorithm
is \emph{memory-less}.

We give an example of how the edge graph of a polytope or,
equivalently, the graph of its normal fan can be
oriented.

\begin{example}
\label{ex: polytope orientation}
Let $P\subset \R^n$ be a polytope whose vertices have positive integer
coordinates and let $\prec$ be a term order on $R$. The following is
an orientation of the edge graph of $P$ without cycles and with a
unique sink: An edge $(p,q)$ is oriented from $p$ to $q$ if and
only if $\x^p\prec\x^q$.
\end{example}

This defines an orientation of the graph of the normal fan of a
polytope for any term order. We would like to mimic this orientation
for any pure full-dimensional fan in $\R^n$. For simplicity we shall
restrict ourselves to fans whose $(n-1)$-dimensional cones allow
rational normals. In view of Propositions \ref{l4} and \ref{l6} this
is no restriction for Gr\"obner fans.

\begin{figure}
\label{fig: not orientable}
\begin{center}
\epsfig{file=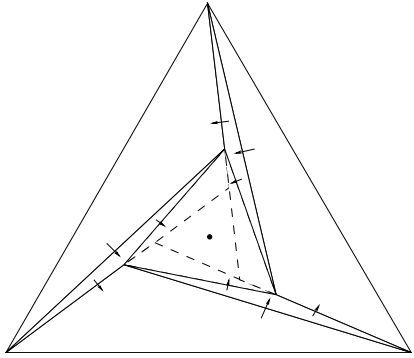,width=5.4cm}
\end{center}
\caption{A fan not having the reverse search property, see Example \ref{ex: not orientable}.}
\end{figure}

\begin{definition}
\label{def:reverse search property}
A pure full-dimensional fan in $\R^n$ is said to have the \emph{reverse search property} if for any
term order $\prec$ the following is an acyclic orientation of its
graph with a unique sink: If $(C_1,C_2)$ is an edge then
$C_1$ and $C_2$ are $n$-dimensional cones with a common facet $F$. Let
$p,q\in\N^n$ such that $q-p\not=0$ is a normal for $F$ with all points
in $C_1\backslash F$ having negative inner product with $q-p$ and all
points in $C_2\backslash F$ having positive inner product with
$q-p$. We orient the edge in direction from $C_1$ to $C_2$ if and only if
$\x^p\prec\x^q$.
\end{definition}
Note that the orientation of an edge in Definition \ref{def:reverse search property} does not depend on the particular choice of $p$ and $q$. Note also that for normal fans of polytopes this orientation agrees with the orientation of the edge graphs of the polytopes in Example \ref{ex: polytope
orientation}.
Not every fan has the reverse search property:
\begin{example}
\label{ex: not orientable}
Figure \ref{fig: not orientable} shows a fan with support $\R_{\geq
0}^3$ intersected with the standard simplex. The intersection is the
non-dotted part of the figure. For every shared $2$-dimensional facet
the orientation of its edge with respect to a term order of form $\prec_{(1,1,1)}$
is indicated by an arrow. The graph has a cycle. The reason is that
the vector $(1,1,1)$ is in the interior of the cone over the dotted triangle and therefore induces the shown orientation with any tie-breaking.
\end{example}

Example \ref{ex: polytope orientation} on the other hand shows that
any normal fan of a polytope has the reverse search property. If $I$ is
a homogeneous ideal the Gr\"obner fan of $I$ is known to be the normal
fan of the \emph{state polytope} of $I$, see \cite{sturmfels} for a
proof. (We should mention that in \cite{mall} it was proven that this is only true if we use the state polytope definition in \cite{sturmfels} and not true with the original definition in \cite{bayer}.) As a consequence the Gr\"obner fan will
have the reverse search property. The reverse search orientation of a
fan with respect to any term order can be carried out on any fan
covering $\R_{\geq 0}^n$ and being the normal fan of a
polyhedron. Since the restricted Gr\"obner fan of any $0$-dimensional
or principal ideal satisfies these conditions it is clear that
these fans have the reverse search property.

It is shown in \cite{jensen} that this line of reasoning cannot be applied to Gr\"obner fans in general. 
In particular, an ideal is presented
whose restricted Gr\"obner fan is not the normal fan of a
polyhedron. For this reason we need a non-trivial argument to prove
the following theorem:
\begin{theorem}
\label{thm:gfan has reverse search property}
The Gr\"obner fan of any ideal $I\subseteq R$ has the reverse search property.
\end{theorem}
The proof is given in the next section.
In Section \ref{section: computation} we will argue that all parts of Algorithm  \ref{alg: reverse search} (finding adjacent edges, finding adjacent vertices and finding search edges) can be implemented efficiently for Gr\"obner fans.

\subsection{Proof: The Gr\"obner fan has the reverse search property}
In this section we prove Theorem \ref{thm:gfan has reverse search property}.
We start by recalling how the polynomial ring can be graded by semigroups. This leads to a more general notion of homogeneous ideals.
\begin{definition}
By a \emph{grading} on $R=k[x_1,\dots,x_n]$ we mean a pair $(A,\A)$
consisting of an abelian semigroup $A$ and a semigroup homomorphism:
$$\A:\N^n\rightarrow A$$ such that $\A^{-1}(a)$ is finite for all $a\in A$. The $\A$-\emph{degree} of a term $c\x^b$ is
$\A(b)$. A polynomial is $\A$-\emph{homogeneous} if all its terms have the
same $\A$-degree. An ideal is $\A$-\emph{homogeneous} if it is generated by a
set of $\A$-homogeneous polynomials.
\end{definition}
For a grading $(A,\A)$ on $R$ we get the direct sum of $k$-vector spaces
$$R=\bigoplus_{a\in A}R_a$$
where $R_a$ denotes the $k$-subspace of $R$ consisting of $\A$-homogeneous polynomials of degree $a$.
Any reduced Gr\"obner basis of an $\A$-homogeneous ideal $I$ consists of $\A$-homogeneous polynomials. In particular, by generalizing the argument of Lemma \ref{l2.2} we get the direct sum
$$I=\bigoplus_{a\in A}I_a$$
where $I_a$ denotes the $k$-subspace of $I$ consisting of $\A$-homogeneous polynomials of degree $a$. The $\A$-homogeneous part $I_a$ is a $k$-subspace of $R_a$.
We define the $\A$-\emph{graded Hilbert function}:
\begin{eqnarray}
H_{I,\A}:A&\rightarrow& \N\\
a&\mapsto& \textup{dim}_k(R_a/I_a)
\end{eqnarray}
\begin{remark}
For a monomial ideal $I$ the standard monomials of degree $a$ form a basis for $R_a/I_a$. Hence $H_{I,\A}(a)$ counts the number of standard monomials of degree $a$.
\end{remark}
In general, as the following well-known proposition shows, the Hilbert function can be found by looking at a monomial initial ideal:
\begin{proposition}
\label{prop:hilbert}
Let $I$ be an $\A$-homogeneous ideal and $\prec$ a term order then
$$H_{I,\A}=H_{\init_\prec(I),\A}.$$
\end{proposition}
\begin{proof}
The linear map taking a polynomial to its unique normal form by the division algorithm on $\G_\prec(I)$ induces an isomorphism of $k$-vector spaces
$$R_a/I_a~\rightarrow~ R_a/\init_\prec(I)_a.$$
\end{proof}

Consider a shared facet of the cones $C_1$ and $C_2$ in the Gr\"obner
fan with a relative interior point $v$. The ``edge ideal''
$\init_v(I)$ is homogeneous with respect to any vector in the relative
interior of the facet and consequently also homogeneous with respect
to any vector in the span of the facet. Since $C_1$ and $C_2$ both
contain positive vectors, so does $\span_\R(C_v(I))$. Recall that
$C_v(I)$ is the closure of the equivalence class of $v$.  Pick a basis
$u_1,\dots,u_{n-1}\in \N^n$ for $\span_\R(C_v(I))$ with $u_1$ being a
positive vector. The vectors induce a grading
$\A_v:\N^{n}\rightarrow\N^{n-1}$ on $R$ by $$\A_v(b)=(\langle
u_1,b\rangle,\dots,\langle u_{n-1},b\rangle)$$
for $b\in\N^n$. The
initial ideal $\init_v(I)$ is $\A_v$-homogeneous.

\begin{lemma}
\label{lem:orientation32}
Let $\prec$ be a term order, $I$ an ideal, $(C_1,C_2)$ a directed edge with respect to the orientation in Definition \ref{def:reverse search property} and $M_1$ and $M_2$ the initial ideals of $C_1$ and $C_2$ respectively. Let $v$ be a relative interior point in the shared facet. Then $\init_\prec(\init_v(I))=M_2$.
\end{lemma}
\begin{proof}
Choose a positive interior point $\omega_2$ of $C_2$. We claim that the
following identities hold:
$$M_2
=\init_{\omega_2}(I)
=\init_{\prec_{\omega_2}}(\init_{\omega_2}(I))
=\init_{\prec_{\omega_2}}(I)
=\init_{\prec_{\omega_2}}(\init_v(I))
=\init_{\prec}(\init_v(I)).$$
The first one holds by the choice of $\omega_2$. The second one is clear since $\init_{\omega_2}(I)$ is a monomial ideal. The third one holds
by Lemma \ref{l2} and Lemma \ref{l6.1}. By Lemma \ref{l2} the fourth equality holds since $v\in C_{\prec_{\omega_2}}(I)=C_{\omega_2}(I)$. To prove the last equality we look at the reduced Gr\"obner basis $\G_\prec(\init_v(I))$. If we can show that $\init_{\prec_{\omega_2}}(g)=\init_\prec(g)$ for
all elements $g\in\G_\prec(\init_v(I))$ then we know that $\G_\prec(\init_v(I))$ is also a Gr\"obner basis with respect to $\prec_{\omega_2}$ and the generators for the initial ideal $\init_{\prec_{\omega_2}}(\init_v(I))$ are exactly the same as those for $\init_\prec(\init_v(I))$. This would complete the proof.

The reduced Gr\"obner basis $\G_\prec(\init_v(I))$ is
$\A_v$-homogeneous. For an element $g$ this implies that the
difference between two of its exponent vectors must be perpendicular
to the shared facet. By Definition \ref{def:reverse search property}
there exists a normal $q-p$ of the facet with $\x^p\prec\x^q$ and
$\langle \omega_2,q-p\rangle>0$. Since $\prec$ and $\prec_{\omega_2}$
agree on one normal vector they must agree on all exponent differences
of elements in $\G_\prec(\init_v(I))$.
\end{proof}

Notice
that by Proposition \ref{prop:hilbert} any initial ideal
$\init_\prec(\init_v(I))$ of $\init_v(I)$ has the same $\A_v$-graded Hilbert
function as $\init_v(I)$.

By a flip we mean a move from one vertex in the graph to
a neighbor. For a degree $a\in\N^{n-1}$ we call $\A^{-1}_v(a)$ the \emph{fiber} over $a$. The $\A_v$-graded Hilbert function of an initial ideal $\init_\prec(\init_v(I))$ counts the number of standard monomials inside each fiber. A flip preserves the Hilbert function. We may think of this as monomials in the monomial initial ideal moving around in the fiber. We wish to keep track of how the
monomials move when we walk in the oriented graph. We define exactly what we mean by ``moving around'':

\begin{definition}
\label{definition phi edge}
Let $\prec, M_1, M_2, u_1, \dots ,u_{n-1}$ and $v$ be as above with \\
$\init_\prec (\init_v (I))=M_2$. Let $N_1$ and $N_2$ be the monomials in $M_1$
and $M_2$ respectively. We define the bijection $\phi_{\prec M_1
M_2}:N_1 \rightarrow N_2$ in the following way: For a monomial
$\x^b\in N_1$ look at the monomials $B_1\subseteq N_1$ and $B_2\subseteq
N_2$ with the same $\A$-degree as
$\x^b$. Since taking initial ideals preserves the
$\A$-graded Hilbert function, $|B_1|=|B_2|$. Sort $B_1$ and $B_2$ with respect to
$\prec$. The bijection $\phi_{\prec M_1 M_2}$ is now defined by taking the first
element of $B_1$ to the first element of $B_2$, the second element of
$B_1$ to the second element of $B_2$ and so on.
\end{definition}

The following lemma is from \cite[Lemma 4.1]{lauritzen}:

\begin{lemma}
\label{lem:lau}
Let $\leq_1$ and $\leq_2$ be two term orders. If $f_1^1,\dots,f_s^1$ is a
vector space basis for $I_a$ such that
$\init_{\leq_1}(f_1^1),\dots,\init_{\leq_1}(f_s^1)$ is a basis for
$\init_{\leq_1}(I)_a$, then there exists a basis $f_1^2,\dots,f_s^2$ for
$I_a$ such that $\init_{\leq_2}(f_1^2),\dots,\init_{\leq_2}(f_s^2)$ is a basis
for $\init_{\leq_2}(I)_a$ and
$$ \init_{\leq_2}(f_1^2)\leq_1 \init_{\leq_1}(f_1^1)$$
$$ \vdots$$
$$ \init_{\leq_2}(f_s^2)\leq_1 \init_{\leq_1}(f_s^1).$$
\end{lemma}

\begin{corollary}
\label{cor:lau}
Let the setting be as in Definition \ref{definition phi
edge}. If $\x^b\in M_1$ then $\phi_{\prec M_1 M_2}(\x^b)\not\prec\x^b$.
\end{corollary}
\begin{proof}
Let $a$ be the $\A$-degree of $\x^b$. We apply Lemma \ref{lem:lau} with $I$ in
the lemma being $\init_v(I)$. Let $\leq_1$ be $\prec$ and $\leq_2$ be
the refinement of the preorder induced by $u_1$ with the reversed order of $\prec$. By the orientation of the graph
$M_1=\init_{\leq_2}(\init_v(I))$ and $M_2=\init_{\leq_1}(\init_v(I))$. By multiplying elements of $\G_\prec(\init_v(I))$ by monomials we can construct a $k$-basis $f^1_1, \dots, f^1_s$ of $\init_v(I)_a$ with
$\init_{\leq_1}(f^1_1),\dots,\init_{\leq_1}(f^1_s)$ being a basis of
$(M_2)_a$. By the lemma there is a basis
$\init_{\leq_2}(f^2_1),\dots,\init_{\leq_2}(f^2_s)$ of $(M_1)_a$. Sort the
list of inequalities in the lemma with $\init_{\leq_2}(f^2_i)$ decreasing
w.r.t. $\prec$ ($\leq_1$). The right hand side can now be sorted with
respect to the same order without violating the inequalities. To see
this use the bubble sort algorithm --- when two adjacent inequalities
are swapped \ldots

\begin{center}
\begin{tabular}{ccccccc}
$ \init_{\leq_2}(f^2_i)$ & $\leq_1$ & $\init_{\leq_1}(f^1_i)$ & & $ \init_{\leq_2}(f^2_{i})$ & $\leq_1$ & $\init_{\leq_1}(f^1_{i+1})$\\
$ \vee_1$ & & $\wedge_1$ & $\mapsto $ & $ \vee_1$ & & $\vee_1$ \\
$ \init_{\leq_2}(f^2_{i+1})$ & $\leq_1$ & $\init_{\leq_1}(f^1_{i+1})$ & & $ \init_{\leq_2}(f^2_{i+1})$ & $\leq_1$ & $\init_{\leq_1}(f^1_{i})$
\end{tabular}
\end{center}
\noindent
\ldots the relations on the right hand side of the arrow hold by transitivity of $\leq_1$.

After sorting, $\x^b$ appears somewhere on the left and $\phi_{\prec
M_1 M_2}(\x^b)$ on the right in the same inequality. This completes the
proof.
\end{proof}

\noindent
{\it Proof of Theorem \ref{thm:gfan has reverse search property}.}
Suppose $C_1,C_2,\dots,C_m$ was a path in the oriented graph with
$C_1=C_m$.  Let $M_1,\dots,M_m$ denote the initial ideals and
$N_1,\dots,N_m$ their monomials. We will prove that the bijection
$\phi:=\phi_{\prec M_{m-1} M_m} \circ \dots \circ \phi_{\prec M_1 M_2}$ is the
identity on $M_1$. Suppose it is not the identity and let $\x^b$ be
the smallest element in $M_1$ with respect to $\prec$ that is not fixed by $\phi$. By Corollary \ref{cor:lau}, $\x^b$
is the image of a smaller element in $M_1$ with respect to $\prec$. But this element is fixed by the minimality of $\x^b$ --- a
contradiction. The composition being the identity implies by Corollary \ref{cor:lau} that $\phi_{\prec M_{i} M_{i+1}}$
is the identity for all $i$. Hence $M_i=M_{i+1}$, contradicting that
$M_1,M_2,\dots,M_m$ is a path.

We claim that $C_\prec(I)$ is the unique sink. If $v$ is in the
relative interior of a facet of $C_\prec(I)$ then by Lemma \ref{l2}
$\init_\prec(\init_v(I))=\init_{\prec}(I)$. By Lemma \ref{lem:orientation32} this means that all edges connected to $C_\prec(I)$ are
ingoing. Hence $C_\prec(I)$ is a sink.

To prove uniqueness let $C_{\prec'}(I)$ be some sink in the oriented
graph. By \cite{MoRo} $\prec$ has a matrix representation
$(\tau_0,\dots,\tau_{n-1})\in\R^{n\times n}$ such that
$\tau_\varepsilon:=\tau_0+\varepsilon
\tau_1+\dots+\varepsilon^{n-1}\tau_{n-1}\in \textup{int } C_\prec(I)$
for $\varepsilon>0$ sufficiently small. Furthermore, for any $f\in R$,
$\init_{\tau_\varepsilon}(f)=\init_\prec(f)$ for $\varepsilon>0$
sufficiently small. If $C_{\prec'}(I)$ is a sink then according to
Definition \ref{def:reverse search property} there exists a complete
list of inner normals $q_1-p_1,\dots,q_r-p_r$ of $C_{\prec'}(I)\cap
\R_{\geq 0}^n$ such that $\init_\prec(\x^{q_i}-\x^{p_i})=\x^{q_i}$. Since
$\tau_\varepsilon$ and $\prec$ pick out the same initial forms on a
finite set of polynomials for $\varepsilon>0$ sufficiently small we
see that $\langle \tau_\varepsilon,q_i\rangle >\langle
\tau_\varepsilon, p_i\rangle$ or, equivalently, $\tau_\varepsilon\in
\textup{int } C_{\prec'}(I)$ for $\varepsilon>0$ sufficiently
small. We conclude that $C_{\prec'}(I)=C_{\prec}(I)$.
$\Box$

\section{Implementation issues}
\label{section: computation}
\label{section: local computation}

We can find a single Gr\"obner cone by applying Buchberger's algorithm and Corollary \ref{l2.1} for some term order.
Since the graph of the Gr\"obner fan of $I$ is connected we may choose any graph traversal algorithm for computing the full dimensional Gr\"obner cones. To do the local computations we need to be able to find the edges (connecting facets) of a full dimensional cone and we need to be able to find the neighbor along an edge. We will see how to do this in the following sections.

Throughout the graph enumeration process we will represent the Gr\"obner cones by their marked reduced Gr\"obner bases, rather than by their defining inequalities, their term orders etc..
This choice is justified by the following known theorem which we shall not prove:
\begin{theorem}
Let $I\subseteq R=k[x_1,\dots,x_n]$ be an ideal. The marked reduced Gr\"obner bases of $I$, the monomial initial ideals of $I$ (w.r.t. a positive vector) and the full-dimensional Gr\"obner cones are in bijection. 
\end{theorem}

An important issue when implementing the algorithms is to identify shared facets. We say that a facet is \emph{flippable} if its relative interior contains a positive vector. The flippable facets in a Gr\"obner fan are always shared. With the right definition of search edges the search tree will only consist of flippable facets.

At the end of the section we will see how the search edge computation in the reverse search algorithm can be implemented and we will explain how to take advantage of symmetry in a Gr\"obner fan traversal.
\subsection{Finding facets}

Suppose that we know a marked reduced Gr\"obner basis $\G_\prec(I)$
with respect to some unknown term order $\prec$. Proposition \ref{l4} (or
Corollary  \ref{l2.1}) tells us how to read off the defining inequality system for
$C_\prec(I)$.

Since $C_\prec(I)$ is full-dimensional the system contains no equations but only inequalities. Some of these inequalities are equivalent in the sense that they are multiples of each other. Taking just one inequality from each equivalence class the problem is now to find irredundant facet normals of a cone --- or equivalently
to find the extreme rays of the dual cone. Checking if a ray is extreme can be done by linear programming.

Not all of the remaining inequalities are guaranteed to define flippable facets. One way to ensure that we only get flippable facets is by adding the constraints $e_i\cdot x\geq 0$ for $i=1,\dots,n$ and ignoring the facets defined by these.

A more efficient method (on some examples) is to find all facets and then remove the non-flippable irredundant facet normals by explicit checks. In our implementation this is done by checking if the inequality system with the inequality in question inverted still has a positive solution. 

As mentioned in \cite{huber} there is an algebraic test that helps us
eliminate redundant inequalities of $C_\prec(I)$. Let $\alpha\in\R^n$
be a coefficient vector of an inequality. If $\alpha$ indeed is
irredundant and defines a facet with a relative interior point $v$
then Corollary \ref{l3} tells us how to compute
$\G_\prec(\init_v(I))$. This marked reduced Gr\"obner basis can be
computed from $\G_\prec(I)$ as $\{\init_v(g)\}_{g\in\G_\prec(I)}$ if
we just know $\alpha$ and not necessarily $v$, see the next section. A
necessary condition for $\alpha$ to be irredundant is that the
computed set $\{\init_v(g)\}_{g\in\G_\prec(I)}$ indeed is a marked
Gr\"obner basis i.e. all S-polynomials reduce to zero. This check even
works for $v$ outside the positive orthant. A quicker necessary
condition that we can check is that every non-zero S-polynomial should
have at least one of its terms in $\init_\prec(I)$. For huge sets of
inequalities the test works extremely well --- 500 inequalities might
reduce to 50 of which maybe 10 are irredundant. Our experience is that
having this test as a preprocessing step can be much faster than
solving the full linear programs with exact arithmetic.

\subsection{Local change}
Let $\G_\prec(I)$ be a known marked Gr\"obner basis and let $F$ be a flippable facet of $C_\prec(I)$. We let $\flip(\G_\prec(I),F)$ denote the unique reduced Gr\"obner basis different from $\G_\prec(I)$ whose Gr\"obner cone also has $F$ as a facet. We will describe an algorithm for computing $\flip(\G_\prec(I),F)$ given $\G_\prec(I)$ and an inner normal vector $\alpha$ for $F$.
For a marked Gr\"obner basis $\G$ and a polynomial $g$ we let $g^\G$ denote the normal form of $g$ modulo $\G$ and note that this form does not depend on the term order but only on $\G$.

\begin{algorithm} $ $\\
{\bf Input:} A marked reduced Gr\"obner basis $\G_\prec(I)$ with $\prec$ being an unknown term order and an inner normal vector $\alpha$ of a flippable facet $F$ of $C_\prec(I)$.\\
{\bf Output:} $\G=\flip(\G_\prec(I),F)$.\\
$\{$\\
\TAB Let $v$ be a positive vector in the relative interior of $F$;\\
\TAB Compute $\G_\prec(\init_v(I))=\{\init_v(g):g\in\G_\prec(I)\}$;\\
\TAB Compute the marked basis $\G_{\prec_{-\alpha}}(\init_v(I))$ from $\G_\prec(\init_v(I))$\\
\TAB \TAB using Buchberger's algorithm;\\
\TAB $\G:=\{g-g^{\G_\prec(I)}:g\in \G_{\prec_{-\alpha}}(\init_v(I))\}$;\\
\TAB Mark the term $\init_{\prec_{-\alpha}}(g)$ in each element $g-g^{\G_\prec(I)}$ in $\G$;\\
\TAB Turn $\G$ into a reduced basis;\\
$\}$
\end{algorithm}
The algorithm is a special case of the local change procedure for a single step in the Gr\"obner walk \cite{collart}. See \cite[Proposition 3.2]{genericgroebnerwalk} for a new treatment and a proof. Here we will just add a few comments on our special case --- the case where $F$ is a facet and not a lower dimensional face:

For any vector $\omega$ in the relative interior of $F$, 
$\init_\omega(I)=\init_v(I)$ is homogeneous with respect to the
$\omega$-grading. Since $F$ is $(n-1)$-dimensional, $\init_v(I)$ is
homogeneous with respect to all vectors inside
$\span_\R(\alpha)^\perp$. All Gr\"obner bases of $\init_v(I)$ are homogeneous in the same
way. Consequently, each of them must consist of polynomials of the form $\sum_{s=0}^t c_s
\x^{(a+sb)}$ where $a\in \N^n$ and $b\in\Z^n$ is parallel to $\alpha$. The
same is true for all polynomials appearing in any run of Buchberger's
algorithm starting from one of these sets. A consequence is that in
order to run Buchberger's algorithm we only need to decide if
we are in the situation where $\x^\gamma\prec \x^{\gamma+\alpha}$ for $\gamma\in\N^n$ or in the situation where $\x^{\gamma+\alpha}\prec \x^\gamma$ for $\gamma\in\N^n$.
Thus specifying $\alpha$ or $-\alpha$ as a term order suffices --- no tie-breaker is needed. The initial ideal $\init_v(I)$ can have at most
two reduced Gr\"obner bases. Both term orders are legal since
$\init_v(I)$ is homogeneous with respect to the strictly positive vector
$v$.  

The Gr\"obner basis $\G_\prec(\init_v(I))$ can be read off from the
marked Gr\"obner basis $\G_\prec(I)$ by taking initial forms of the
polynomials with respect to $v$, see Corollary \ref{l3}. Taking the initial form $\init_v(g)$ of a polynomial $g\in\G_\prec(I)$ without computing $v$ is done as follows. By Corollary \ref{l2.1}, $\init_\prec(\init_v(g))=\init_{\prec}(g)$ and thus we already know one term of $\init_v(g)$ since $\init_\prec(g)$ is the marked term of $g$ in $\G_\prec(I)$. Since every $\omega$ in the relative interior of $F$ will have $\G_\prec(\init_\omega(I))=\G_\prec(\init_v(I))$ the remaining terms of $\init_v(g)$ are exactly the terms in $g$ with the
same $\omega$-degree as $\init_\prec(g)$ for all $\omega$ in the
relative interior of $F$ and consequently for all $\omega$ in
$\span_\R(\alpha)^\perp$. In other words a term of $g$ is in $\init_v(g)$ if and only if its exponent vector minus the exponent of $\init_\prec(g)$
is parallel to $\alpha$. The term order $\prec$ does not
have to be known for this step, nor does it have to be known in the computation
of $\G_{\prec_{-\alpha}}(\init_v(I))$ or in any other subsequent
step. The vector $v$ also remains unknown in the entire process.

\subsection{Computing the search edge}
\label{subsec: search edge}
Let $\prec$ be the term order used for orienting the graph of the Gr\"obner fan. In Algorithm \ref{alg: reverse search} the search edge $e_{C_{\prec'}(I)}$ has to be computed given $\G_{\prec'}(I)$ where ${\prec'}$ is some unspecified term order. According to Proposition \ref{proposition: reverse search graph} the definition of search edges can be arbitrary.
However, efficiently computing a search edge requires a good definition. Our search edges will always come from flippable facets.

One strategy for locally computing the search edge $e_{C_{\prec'}(I)}$ is to compute a unique representation
of each flippable facet of the Gr\"obner cone $C_{\prec'}(I)$ and then choose the smallest of
these facets to be $e_{\C_{\prec'}(I)}$ in some lexicographic order. This method
requires all facets to be computed every time we check if
``$e_{u}=(u,v)$'' in Algorithm \ref{alg: reverse search}.

A better strategy is to draw a straight line from a point in the cone
$C_{\prec'}(I)$ to the cone of the sink and choose the first facet intersecting
this line as $e_{C_{\prec'}(I)}$.  A point in the cone $C_{\prec'}(I)$ can be computed
deterministically by linear programming. Two problems arise. The
straight line might not intersect a unique facet and we may not know a
point in the cone of the sink. Both problems can be solved using
formal perturbation of the end points of the line. This was worked
out in detail in \cite{genericgroebnerwalk}. Here we explain
how it works for lexicographic term orders and with one end point perturbed.

\begin{lemma}
Let $I\subseteq R$ be an ideal and $\prec$ the lexicographic term
order with $x_1\succ x_2\succ \dots \succ x_n$. Define $\tau_\varepsilon=(\varepsilon^0,\varepsilon^1,\dots,\varepsilon^{n-1})$. There exists a
$\delta>0$ such that
$\init_{\tau_\varepsilon}(I)=\init_\prec(I)$
for all $\varepsilon\in (0,\delta)$.
\end{lemma}
\begin{proof}
This follows from Lemma \ref{l1} since $\prec$ and
$\tau_\varepsilon$ agree on a finite set of polynomials for small
$\varepsilon$.
\end{proof}

Let $\sigma$ be a deterministically computed interior point of the
cone of $C_{\prec'}(I)$ and assume for simplicity that $\sigma\in\N^n$. For
sufficiently small $\varepsilon>0$ the line segment
$$\omega(t):=(1-t)\sigma+t\tau_\varepsilon \textup{ with }t\in[0,1]$$
intersects a facet of $C_{\prec'}(I)$ unless $C_{\prec'}(I)$ is the sink.

Let $\{\alpha_1,\dots,\alpha_m\}$ be the set of potential inner facet normals read off
from $\G_{\prec'}(I)$.  We
are only interested in the vectors $\alpha_i$ where
$\langle\sigma,\alpha_i\rangle >0$ and $\langle\tau_\varepsilon,\alpha_i\rangle
<0$.  Let $t_i$ denote the $t$-value for the intersection of the line
segment and the hyperplane defined by $\alpha_i$. Then
$$ t_i:={\langle \sigma, \alpha_i\rangle\over\langle\sigma,\alpha_i\rangle -\langle\tau_\varepsilon,\alpha_i\rangle}.$$
We wish to find $i$ such that $t_i$ is smallest (for small $\varepsilon$).
\begin{eqnarray}
t_i&<&t_j ~\Longleftrightarrow \\
 {\langle \sigma, \alpha_i\rangle\over\langle\sigma,\alpha_i\rangle -\langle\tau_\varepsilon,\alpha_i\rangle} &<& {\langle \sigma, \alpha_j\rangle\over\langle\sigma,\alpha_j\rangle -\langle\tau_\varepsilon,\alpha_j\rangle} ~\Longleftrightarrow \\
 {\langle\sigma,\alpha_i\rangle -\langle\tau_\varepsilon,\alpha_i\rangle\over \langle \sigma, \alpha_i\rangle} &>& {\langle\sigma,\alpha_j\rangle -\langle\tau_\varepsilon,\alpha_j\rangle\over \langle \sigma, \alpha_j\rangle} ~\Longleftrightarrow \\
{\langle\tau_\varepsilon,\alpha_i\rangle\over \langle \sigma, \alpha_i\rangle} &<& {\langle\tau_\varepsilon,\alpha_j\rangle\over \langle \sigma, \alpha_j\rangle} ~\Longleftrightarrow \\
\langle\tau_\varepsilon,{\langle \sigma, \alpha_j\rangle \alpha_i}\rangle &<& \langle\tau_\varepsilon,{\langle \sigma, \alpha_i\rangle\alpha_j}\rangle ~\Longleftrightarrow \\
\x^{\langle \sigma, \alpha_j\rangle \alpha_i} &\prec& \x^{\langle \sigma, \alpha_i\rangle \alpha_j}
\end{eqnarray}
We see that for $\varepsilon$ sufficiently small ``$t_i<t_j$'' does not depend on $\varepsilon$.
Furthermore, there cannot be any ties, unless $\alpha_i$ and $\alpha_j$ represent the same hyperplane.
This gives an easy method for defining and computing $e_{C_{\prec'}(I)}$.
We simply choose the facet defined by $a_i$ where $t_i$ is smallest among $\{t_1,\dots,t_m\}$ (for small $\varepsilon>0$).

\subsection{Exploiting symmetry}
\label{section: symmetry}
In this section we explain how to take advantage of symmetry to
speed up computations.
The symmetric group $S_n$ acts on
polynomials and ideals of $R$ by permuting variables
and on $\R^n$ by permuting coordinate entries. Let $I\subseteq R$ be an
ideal. We call a subgroup $\Gamma\leq S_n$ a \emph{symmetry group} for $I$
if $\pi(I)=I$ for all $\pi\in \Gamma$.  If we know a symmetry group
for $I$ we can enumerate the reduced Gr\"obner bases of $I$ up to
symmetry.
Let $\Gamma$ be such a symmetry group for $I$.

In our description all Gr\"obner bases will be marked
and reduced. Thereby each one will uniquely represent its initial
ideal and Gr\"obner cone. For a Gr\"obner basis $\G$ of $I$ we use the notation
$\Gamma_\G=\{\pi(\G)\}_{\pi\in \Gamma}$ for its orbit.

The idea is to exploit the identity
$\flip(\pi(\G),\pi(F))=\pi(\flip(\G,F))$ for all $\pi\in\Gamma$. In other words $\Gamma$ is a group of automorphisms of the graph of the Gr\"obner fan of $I$. The \emph{quotient graph} is defined to be the graph whose vertices are the orbits of Gr\"obner bases with two orbits 
$\Gamma_\G$ and $\Gamma_{\G'}$ being connected if there exists a facet $F$ of the Gr\"obner cone of $\G$ such that $\flip(\G,F)\in \Gamma_{\G'}$.
The flip graph may have loops.

The symmetry-exploiting algorithm enumerates the quotient graph by a
breadth-first traversal.  Orbits are represented by Gr\"obner basis
representatives. One question that arises is how to check if two
Gr\"obner bases $\G$ and $\G'$ represent the same orbit. A solution is
to run through all elements $\pi\in \Gamma$ and check if $\pi(\G)$
equals $\G'$, or even better to make a similar check for the monomial
initial ideals.  Although this does not seem efficient, it is still
much faster in practice than redoing symmetric Gr\"obner basis and
polyhedral computation as we have done in the usual reverse search or
breadth-first enumeration without symmetry. It is not clear how to
combine symmetry-exploiting with reverse search.

\def\FACETS{\ensuremath{{T_{\bf{facets}}}}}
\def\SHOOT{\ensuremath{{T_{\bf{shoot}}}}}
\def\FLIP{\ensuremath{{T_{\bf{flip}}}}}
\def\LP{\ensuremath{{T_{\bf{lp}}}}}

\section{Complexity}
\label{section: complexity}
In this section we will discuss the complexity of enumerating the maximal cones of the Gr\"obner fan of an ideal $I$ by reverse search. We will assume that $I$ is homogeneous with respect to a positive vector. This guarantees that any facet of a full-dimensional Gr\"obner cone is flippable.

We identify the following important sub-algorithms:
\begin{itemize}
\item Computation of the facet normals of the Gr\"obner cone of a marked reduced Gr\"obner basis $\G$. We denote the time for this operation by \FACETS($\G$).
\item Computation of a search edge given a marked reduced Gr\"obner basis $\G$ as described in Subsection \ref{subsec: search edge}. We denote the time for this operation by \SHOOT($\G$).
\item Conversion of a marked reduced Gr\"obner basis $\G_1$ into a marked reduced Gr\"obner basis $\G_2$ where the Gr\"obner cones of $\G_1$ and $\G_2$ are assumed to share a facet. We denote the time for this operation by \FLIP($\G_1$,$\G_2$).
\end{itemize}
For simplicity we will assume that the time used for solving any of the linear programs in our algorithms only depends on the dimensions of its matrix form. We let \LP($n$,$r$) be the time needed to solve a linear programming problem with $n$ variables and $r$ constraints.
Then \FACETS($\G$) and \SHOOT($\G$) can be expressed in terms of \LP($n$,$r$).

The time \FACETS($\G$)$\in O(\LP(n,r)r)$ where $r$ is the number of non-leading terms in $\G$. The reason is that each non-leading term in $\G$ gives an inequality in the description of the Gr\"obner cone. Checking if the inequality defines a facet takes one linear program. In addition duplicates should be removed from the set of facet normals and further vectors should be eliminated until no parallel vectors exist. The time for this step is dominated by the time for solving LPs. 

The time \SHOOT($\G$)$\in O(r n^2+\LP(n,r))$ where $r$ is the number of non-leading terms in $\G$. The first step in the algorithm is to deterministically find a relative interior point of the Gr\"obner cone. This is done in time $\LP(n,r)$. After this the smallest vector among the $r$ defining vectors for the cone with respect to the ordering in Subsection \ref{subsec: search edge} needs to be found. Comparing two vectors takes $O(n^2)$ operations in the worst case. These are operations in $\Q$. In the above estimate we assume that each operation takes constant time. 

We have no good bound for the complexity of flipping.
Now we count the number of times each of the three sub-algorithms are applied when enumerating the graph of the Gr\"obner fan of $I$ using reverse search.

\begin{itemize}
\item The facets of each Gr\"obner cone are computed exactly once in Algorithm \ref{alg: reverse search} (right after the Gr\"obner basis has been output). We remark that since we are only interested in facets with the correct orientation the number of LPs that really need to be solved is lower than the $r$ in the discussion above. We will not take this into account in our analysis.
\item Checking if an edge is a search edge is done once for every edge. Every time we need to recompute a search edge and compare it to the edge. Hence the total number of times we need to compute a search edge is equal to the number of edges in the graph of the Gr\"obner fan of $I$.
\item When a vertex $\G$ is processed by Algorithm \ref{alg: reverse search} we must test for every ingoing edge if the edge is a search edge. To test this we first compute $\flip(\G,F)$ where $F$ is the facet of the cone corresponding to the edge in question. If the edge is equal to the search edge of $\flip(\G,F)$ we do an enumeration of the subtree with root $\flip(\G,F)$. If not, $\flip(\G,F)$ is forgotten. Since all vertices are processed once and every edge is ingoing for exactly one vertex the number of times $\flip(\G,F)$ needs to be computed is equal to the number of edges in the graph. We remark that the variant of the reverse search where the search path for the current vertex is not stored on the recursion stack would require twice as many computations of this kind.
\end{itemize}
Let $E$ be the edges and $V$ be the vertices of the graph. The total time complexity of the enumeration of $(V,E)$ is:
$$O(\sum_{\G\in V} \FACETS(\G) +  \sum_{(\G_1,\G_2)\in E} \SHOOT(\G_1) + \sum_{(\G_1,\G_2)\in E}\FLIP(\G_1,\G_2))$$
Substituting with the time needed for solving the LPs we get the following theorem:
\begin{theorem}
Let $(V,E)$ be the graph of the Gr\"obner fan of $I$. The time complexity for computing this graph given a marked reduced Gr\"obner basis is in the class of functions
$$O(\sum_{\G\in V} \LP(n,r(\G))r(\G) +\sum_{(\G_1,\G_2)\in E} r(\G_1) n^2+\LP(n,\G_1) + \sum_{(\G_2,\G_1)\in E}\FLIP(\G_1,\G_2))$$
where $r(\G)$ is the number non-leading terms in the marked reduced Gr\"obner basis $\G$. In particular, the first two terms are bounded by a polynomial in the size of the output.
\end{theorem}

\begin{corollary}
If for a given class of ideals the time $\FLIP(\G_1,\G_2)$ is bounded by a polynomial in the size of the binary encoding of $\G_1$ and $\G_2$ then the enumeration of the reduced Gr\"obner bases for an ideal in the class by reverse search is a polynomial time algorithm in the size of the output. 
\end{corollary}

\section{Computational results and examples}
\label{sec experiments}
The algorithms presented in this paper were implemented in the
software package Gfan \cite{gfan}. In this section we present examples
of Gr\"obner fans computed using this package. The first example comes
with a picture and gives an idea of the kind of geometric shape a
Gr\"obner fan might have.
\begin{figure}
\begin{center}
\epsfig{file=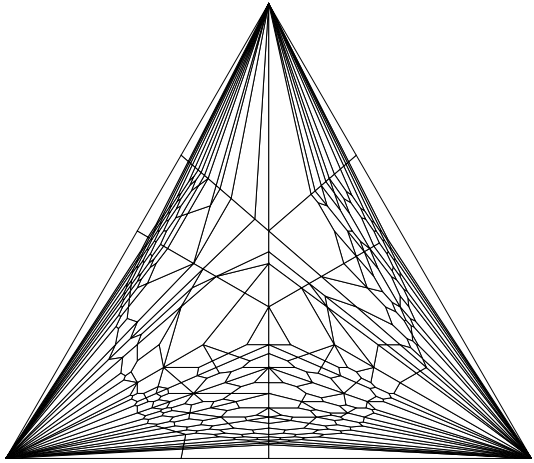,width=7cm}
\end{center}
\caption{The Gr\"obner fan of the ideal in Example
  \ref{ex:sturmfels3.9} intersected with the standard $2$-simplex. The
  $a$-axis is on the right, the $b$-axis on the left and the $c$-axis
  at the top.} 
\label{fig:sturmfels3.9}
\end{figure}

\begin{example}
\label{ex:sturmfels3.9}
\cite[Example 3.9]{sturmfels} Consider the ideal $ I=\langle
a^{5}-1+c^{2}+b^{3}, b^{2}-1+c+a^{2},
c^{3}-1+b^{5}+a^{6}\rangle\subseteq\Q[a,b,c]$. The Gr\"obner fan of $I$
has 360 full-dimensional cones and the Gr\"obner region is $\R_{\geq
  0}^3$. The intersection of the fan with the standard simplex in
$\R^3$ is shown in Figure \ref{fig:sturmfels3.9}.
\end{example}

We now list some families of ideals used in our computations. The
Gr\"obner fans of these ideals have been computed for the parameters
listed in the table of Figure \ref{fig:example table}. The ambient
field is always $\Q$.  The columns of the table are to be interpreted
as follows. In each row, the first column contains the name of the
ideal (to be explained below). The second column lists $n$, the number
of variables in the ideal. The third column lists $h$, the dimension
of the lowest dimensional Gr\"obner cone $C_0(I)$. Note that $h$ is
the dimension of the {\em homogeneity space} of the ideal which is the
common subspace contained in every Gr\"obner cone of the ideal.  The
quantity ``d'' is the lowest total degree of any reduced Gr\"obner
basis of the ideal and ``D'' is the highest.  The $f$-{\em vector}
of the Gr\"obner fan is an ordered list of the number of
$h$-dimensional cones, $h+1$-dimensional cones etc., up to the number
of $n$-dimensional cones.

\begin{figure}
\begin{center}
\begin{tabular}{|l|r|r|r|r|l|}
\hline
Example & $n$ & $h$ & d & D & $f$-vector \\
\hline
$\textup{Det}_{3,3,4}$   & 12 & 6 & 3 & 3 & (1,12,66,204,342,288,96)\\
$\textup{Det}_{3,3,5}$   & 15 & 7 & 3 & 3 & (1,45,585,3390,10710,19890,21750,12960,3240)\\
$\textup{Det}_{3,4,4}$   & 16 & 7 & 3 & 5 & (1,?,?,?,?,?,?,?,?,163032)\\
$\textup{Detsym}_{3,4}$  & 10 & 4 & 3 & 8 & (1,518,5412,20505,36024,29808,9395)\\
$\textup{Grass}_{2,5}$   & 10 & 5 & 2 & 3 & (1,20,120,300,330,132)\\
$\textup{Cyclic}_{5}$    &  5 & 0 & 8 & 15 & (1,?,?,?,?,55320)\\
$\textup{J}_{4}$         &  4 & 1 & 3 & 8 & (1,200,516,318)\\
\hline
\end{tabular}
\end{center}
\caption{Statistics for the Gr\"obner fans computed using Gfan.
}
\label{fig:example table}
\end{figure}

\begin{example}
  Let $\textup{Det}_{t,m,n}$ denote the ideal in the polynomial ring
  in $mn$ variables generated by the $t\times t$ minors of the matrix:
$$ \left( \begin{array}{cccc}
x_{11} & x_{12} & \cdots & x_{1n} \\
x_{21} & x_{22} & \cdots & x_{2n} \\
\vdots & \vdots & \ddots & \vdots \\
x_{m1} & x_{m2} & \cdots & x_{mn} \\
\end{array} \right).$$
\end{example}

\begin{example}
  Let $\textup{Grass}_{d,n}$ denote the ideal in the polynomial ring
  in $n\choose d$ variables generated by the relations on the $d\times
  d$ minors of a $d\times n$ matrix.
\end{example}

\begin{example}
  Let $\textup{Detsym}_{t,n}$ denote the ideal in the polynomial ring
  in $n(n+1)\over 2$ variables generated by the $t\times t$ minors of
  a symmetric matrix of variables. For example,
  $\textup{Detsym}_{3,4}$ is generated by the $3 \times 3$ minors of the
  following matrix:
$$ \left( \begin{array}{cccc}
a & b & c & d \\
b & e & f & g \\
c & f & h & i \\
d & g & i & j \\
\end{array} \right).$$

\end{example}

\begin{example}
  Let $\textup{Cyclic}_5$ denote the ideal $\langle a+b+c+d+e,
  ab+bc+cd+de+ae, abc+bcd+cde+dea+eab, abcd+abce+abde+acde+bcde,
  abcde-1\rangle\subseteq k[a,b,c,d,e].$ In general, $\textup{Cyclic}_n$
  stands for the generalization of this polynomial system to $n$
  variables \cite{BF}. These polynomial systems have become benchmarks
  for computer algebra packages and their lexicographic Gr\"obner bases
  are notoriously hard to compute.
\end{example}

\begin{example}
  Let $K_n$ denote the complete graph on $n$ vertices and $I_n$ be the
  {\em Stanley-Reisner} ideal of this graph. The Stanley-Reisner ideal
  of a simplicial complex $\Delta$ is the ideal generated by all
  monomials $x_{i_1}x_{i_2} \cdots x_{i_n}$ such that $\{i_1, \ldots,
  i_n\}$ is not a face of $\Delta$. Apply a generic linear change of
  coordinates to $I_n$ to obtain the ideal $J_n$. The generators of
  $J_n$ typically have very complicated coefficients. For example, the
  first generator in our $J_4$ was 
$$\begin{array}{rrrrr}
\vspace{0.1cm} 
  a^3&+  \frac{4980248985}{343338664}a^2c&+\frac{2079196217}{257503998}abc&+\frac{86858380}{128751999}b^2c&-\frac{2205648949}{42917333}ac^2\\
\vspace{0.1cm} 
 &-\frac{359584197}{171669332}bc^2 &-\frac{84523033581}{1373354656}c^3&-\frac{11737327991}{51500799600}a^2d&- \frac{16299027451}{38625599700}abd\\
\vspace{0.1cm} 
&+\frac{1194144014}{9656399925}b^2d&+\frac{394500908221}{25750399800}acd &-\frac{47953955497}{25750399800}bcd&+\frac{195491595943}{2985553600}c^2d\\
&-\frac{4583330213}{3862559970}ad^2&+  \frac{181743499}{364392450}bd^2 &-\frac{429736138279}{25750399800}cd^2&+  \frac{8566043731}{12875199900}d^3.\\
\end{array}
$$
The initial ideals of $J_n$ are known as the {\em generic initial
  ideals} of $I_n$. The reverse lexicographic generic initial ideals
of an ideal have played an important role in commutative algebra and
algebraic geometry while other generic initial ideals have not been
explored too much. We computed the Gr\"obner fan of $J_4$.
\end{example}

Extracting the f-vector from the full-dimensional Gr\"obner cones
produced in the enumeration process was the most time-consuming part
of the computation of these examples. In example
$\textup{Det}_{3,3,4}$ this extraction was not possible to complete
within reasonable time with the current software package. For this
particular example the $163032$ full-dimensional Gr\"obner cones were
computed up to the action of a symmetry group of order $576$. The
full-dimensional cones come in $289$ orbits. The computation of the
full dimensional cones up to symmetry took 7 minutes on a 2.4 GHz
Pentium processor. Using reverse search without symmetry the same
computation would take approximately 14 hours. The f-vector extraction
routine in Gfan only works for complete fans. This is why the f-vector
for the $\textup{Cyclic}_5$ example is not shown.

\appendix

\bibliographystyle {hplain}
\bibliography{jensen.bib}

\bigskip
 
\noindent Komei Fukuda, Institute for Operations Research and Institute of Theoretical Computer Science
ETH Zentrum, CH-8092 Zurich, Switzerland
and \break
Mathematics Institute / ROSO
EPFL, CH-1015 Lausanne, Switzerland\break
{\tt fukuda@ifor.math.ethz.ch}.

\bigskip

\noindent Anders N. Jensen, Institut for Matematiske Fag,
Aarhus Universitet, \break DK-8000 \AA rhus, Denmark,
{\tt ajensen@imf.au.dk}.
  
\bigskip
 
\noindent Rekha R. Thomas, Department of Mathematics,
University of Washington, \break Seattle, WA 98195-4350, USA,
{\tt thomas@math.washington.edu}.

\end{document}